# Stabilizing Gain Selection of Networked Variable Gain Controller to Maximize Robustness Using Particle Swarm Optimization

Indranil Pan, Saptarshi Das, Soumyajit Ghosh, Amitava Gupta

*Abstract*—Networked Control Systems (NCSs) are often associated with problems like random data losses which might lead to system instability. This paper proposes a method based on the use of variable controller gains to achieve maximum parametric robustness of the plant controlled over a network. Stability using variable controller gains under data loss conditions is analyzed using a suitable Linear Matrix Inequality (LMI) formulation. Also, a Particle Swarm Optimization (PSO) based technique is used to maximize parametric robustness of the plant.

## I. INTRODUCTION

REAL time control loops which are closed over shared communication networks have received increased attention from researchers due to the significant advantages they offer over traditional control systems [1]. Some of the key features of networked control systems include reduced wiring due to the usage of shared medium as opposed to point to point connection, low cost and modularity, flexibility of architecture etc. NCSs have wide applicability in factory automation, aircraft and space applications, automobiles, remote diagnostics, manufacturing plant monitoring and access in hazardous environments. However with the introduction of the network into the control loop, various practical issues and constraints comes into play which affects the control performance [2]. Among all the network constraints, level of packet drop-out (due to buffer overflows etc.) and random network delays due to transmission, routing etc. introduced within the network plays pivotal role so far as the control loop performance is concerned. Random packet drop-outs and delays can have varying effect on the control performance ranging from degrading the closed loop response to making the system unstable altogether [3].

Theoretical analysis for the stabilization of linear systems over networks in the presence of arbitrary and Markovian packet losses have been studied by [4]. Other studies [5] have analyzed the network as a switched system using the concept of lifted sampling period and derived stability conditions for delays and packet dropouts. Model Predictive control over networks has been studied by [6] and has been shown to work well over the classical state feedback controller for NCS. The predictive control methodology [7]-[12] exploits the fact that for data networks like Ethernet, a single packet of data can accommodate many control signal values and hence predicted values of the future control signals can be clubbed in the present packet and transmitted to the actuator which decides the appropriate control signal to be used. In [13] active compensation for network delay and drop has been done by this method and has shown improved performance over other conventional methods. However, the controller performs an online finite horizon optimization to predict the control signal values at each time instant. This poses a constraint on the real time requirement that the controller must finish the computation of every packet, within a very small fraction of the sampling time.

Instead of the finite horizon optimization at each time instant, this paper proposes static state-feedback gains for each of the future predicted control signals. This reduces the controller complexity and thus the hard limit of the real time requirements can easily be met. The predicted control signal values for the finite horizon are put in the same packet as the present control signal and sent to the actuator buffer. The actuator updates the plant with values from the buffer based on different scenarios of arrival of the packets at the actuator.

In this paper, the NCS is modelled as a set of switched Linear Time Invariant (LTI) systems with uncertain switching between them caused by network conditions like data loss and variable delay. While stability of the NCS with network considerations is ascertained by way of finding a feasible solution for the LMIs using multiple Lyapunov functions, the gains are so chosen such that the parametric robustness of the nominal system is maximum. The nominal system here implies the closed loop system without network considerations (e.g. variable delay and data loss) and the corresponding controller gains are termed as the nominal gains. Other robust stability approaches in NCS can be found in [14]-[15]. The present approach fits the largest n-dimensional hypersphere in an n-dimensional stable controller space and selecting the nominal controller gains corresponding the center of the hypersphere. Particle Swarm Optimization is used to obtain the set of predictive controller gains that maximizes the volume of the hypersphere. This ensures maximum variation in controller parameters and since any uncertainty in the system parameters can be reflected as change in controller gains [16]-[17], the designed NCS allows for maximum system parameter variation due to modeling uncertainty.

The rest of the paper is organized as follows. Section II

Manuscript received April 14, 2011. This work has been supported by the Board of Research in Nuclear Sciences (BRNS) of the Department of Atomic Energy (DAE), India, sanction no. 2009/36/62-BRNS, dated November 2009.

I. Pan, S. Das, S. Ghosh and A. Gupta are with the Department of Power Engineering, Jadavpur University, Kolkata, India. (E-mail: indranil.jj@student.iitd.ac.in, indranil@pe.jusl.ac.in).

deals with the theoretical formulation. Section III presents the main results with a proof of the stability conditions in the Lyapunov sense. Section IV supports the above formulation with an illustrative example. The paper ends in Section V with the conclusions and is followed by the references.

## II. THEORETICAL FORMULATION

The continuous time plant is discretized with a specified sampling time and the corresponding linear discrete time system is given by

$$x(k+1) = Fx(k) + Gu(k) \qquad (1)$$

where, $k \in \mathbb{Z}_+$ is the time step, $x(k) \in \mathbb{R}^n$ is the system state and $u(k) \in \mathbb{R}^m$ is the control input. $x_0 := x(0)$ is the initial state and $F$ and $G$ are two constant matrices of appropriate dimensions. Fig. 1 shows the schematic diagram of the system. The sensor is time driven and samples the plant-states at fixed instants of sampling time. The sensor sends the packets to the controller along with a time stamp. It is assumed that full state information of the plant is available. The controller is a full state feedback controller and is also time driven. Generally the convention is to use an event driven controller which calculates the control signal as soon as it is received and sends it over the network to the actuator [18]. This eliminates the additional time till the next sampling instant that the controller must wait before it transmits the control signal as in the time driven case. But in safety critical applications event driven controllers are not used since specific cases may occur where after a long delay, many control signals come together within the same sampling time thus increasing the instantaneous network load to a high value which is undesirable. Hence the choice of this type of controller is made in the present paper. The predictive state feedback controller is implemented as a set of static gains and the controller embeds the values of the present control signal and the next predicted control signals in the same packet. The buffer stores the received packet and applies the appropriate value of the control signal to the actuator based on some logic. The control signal is given by $u = -K_z x$ where $z \in \{1, 2, ..., M_{drop}\}$ and is applied to the plant by the actuator on a logic based on the arrival of the preceeding packets as discussed in Definition 2. There is a skew between the time line of the sensor, controller and actuator to ensure that the packet sent by the sensor for a time step is picked up by the controller and similarly the packet sent by the controller for a particular time step is picked up by the actuator. Thus if $\tau_1$ be the time skew between the sensor and the controller, $\tau_2$ be the time skew between the controller and the actuator, then the following is assumed to hold:

$T_s > \tau_1 + \tau_2$, where $T_s$ is the sampling time of the system. Thus the timelines have an offset which is assumed to be maintained with a precision so that the controller picks up the current sampled data sent by the sensor and the actuator picks up the current control signal sent by the controller.

In a practical scenario, a number of control loops may use the same medium and the nature of the traffic might be burst periodic in certain networks like Ethernet. This might lead to a jitter in transmissions and variable delay. However a packet which belongs to a particular time step is considered to be dropped if it does not reach the actuator in the same time step. Any out of order packet is discarded by the buffer and considered to be dropped.

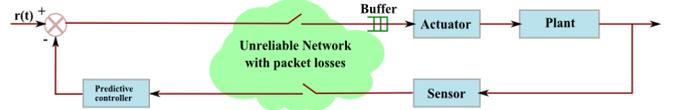

Fig. 1. Schematic of the networked predictive control system.

Let $\omega := \{i_1, i_2, ...\}$ be a subsequence of $\mathbb{Z}_+$ indicating the sequence of time points where there is successful data transmission from the sampler to the actuator buffer. Also let $M_{drop} := \sup_{i_k \in \omega} |i_{k+1} - i_k|$ be the maximum upper bound of consecutive packet losses.

### A. Definition 1

The mathematical model of the packet loss process is defined as

$$\xi(i_k) := \{i_{k+1} - i_k \mid i_k \in \omega\} \qquad (2)$$

which takes values in the finite state space $\mu := \{1, 2, ..., M_{drop}\}$. The packet loss process is assumed to take random values in $\mu$ without following any pre-specified distribution e.g. Markovian, Bernoulli etc.

### B. Definition 2:

The buffer logic which decides the control signal $u$ to the actuator at each sampling instant $k$, is defined as

$$u(k) := \psi_{k-\rho+1}(\rho) \text{ where } \rho = k - i_k + 1.$$

Here, $\psi$ denotes the buffer, the subscript $k - \rho + 1$ denotes the time index, the buffer was last updated by an effective packet. $\rho$ denotes the array index of the buffer from which the control signal is provided to the actuator.

Assuming, $M_{drop} = 3$, i.e. maximum 2 samples can be dropped consecutively, Fig. 2 shows the arrival of the packets at the buffer. The crosses indicate the instants of time where the packets do not reach the actuator buffer. The blue lines indicate the arrival of the packets at those particular instants. The buffer state is shown in Table I. At each sample the control signal in the "Present" column is used to update the plant. Note that the number of control signals that each packet brings is equal to $M_{drop}$. The control signal in each control packet is represented by $u_{pq}$ where $p \in \{1, 2, 3, ...\}$ denotes the sample time at which the packet arrived and $q \in \{1, 2, 3 ..., M_{drop}\}$ represents the predicted values of the next $M_{drop}$ number of control signals from the present sampling instant.

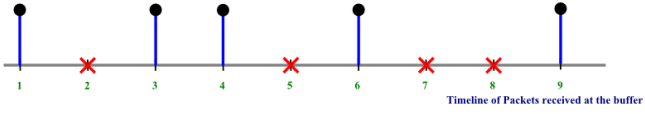

Fig. 2. Schematic showing the timeline of the packets received at the buffer for an example case.

TABLE I
BUFFER STATES CORRESPONDING TO PACKET ARRIVAL IN FIG. 2

| Sample no. | Present | 1st Predicted Sample | 2nd Predicted Sample |
|---|---|---|---|
| 1 | $u_{11}$ | $u_{12}$ | $u_{13}$ |
| 2 | $u_{12}$ | $u_{13}$ | - |
| 3 | $u_{31}$ | $u_{32}$ | $u_{33}$ |
| 4 | $u_{41}$ | $u_{42}$ | $u_{43}$ |
| 5 | $u_{42}$ | $u_{43}$ | - |
| 6 | $u_{61}$ | $u_{62}$ | $u_{63}$ |
| 7 | $u_{62}$ | $u_{63}$ | - |
| 8 | $u_{63}$ | - | - |
| 9 | $u_{91}$ | $u_{92}$ | $u_{93}$ |

## III. DISCRETE TIME SWITCHED SYSTEM STABILITY OF NCS USING PREDICTIVE CONTROLLER

The plant (1) at each instant is updated by the control values in the buffer depending on $\xi(i_k)$. Thus at any given instant the plant can be controlled by a packet index $\psi \in \mu$ in the buffer. Let us consider the augmented state $\Gamma(k) = [x(k)\ x(k-1)\ ...\ x(k-M_{drop}-1)]^T$. The NCS with the state feedback predictive controllers can be cast in the form of a discrete time switched system given by:

$$\Gamma(k+1) = \Phi_{\sigma(k)}\Gamma(k) \quad (3)$$

where, $\sigma(k)$ is a piecewise constant function, known as a switched signal, which takes values in the finite set $\Lambda := \{1, 2, ..., M_{drop}\}$. $\Phi_{\sigma(k)} \in \mathbb{R}^{M_{drop} \times M_{drop}}$ is of the generalized form:

$$\Phi_{\sigma(k)} = \begin{bmatrix} A & \cdots & -BK_{\sigma(k)} & \cdots & 0 \\ I & 0 & \cdots & \cdots & 0 \\ 0 & I & 0 & \vdots & \vdots \\ \vdots & \cdots & \ddots & \vdots & \vdots \\ 0 & \cdots & 0 & I & 0 \end{bmatrix} \quad (4)$$

where,

$$\Phi_{\sigma(k)}(1, \sigma(k)) = \begin{cases} A - BK_1 & \text{for } \sigma(k) = 1 \\ -BK_{\sigma(k)} & \forall \sigma(k) \in \{2, \cdots, M_{drop}\} \end{cases}$$

### A. Theorem 1

The NCS defined by the switched system (3) is asymptotically stable for the arbitrary packet loss process $\xi(i_k)$, if for $i, j \in \Lambda$, there exist positive definite matrices, $P_i$ and $P_j$ satisfying the following set of Linear Matrix Inequalities (LMIs).

$$\Phi_i^T P_j \Phi_i - P_i < 0 \quad (5)$$

where, $\Phi_i$ is of the form (4) with specified controller gains $K_l\ \forall\ l \in \Lambda$.

*Proof:* For the switched system (3) let us define multiple Lyapunov functions of the form (6), for each switched state.

$$V(k) = \Gamma^T(k) P_{\sigma(k)} \Gamma(k) \quad (6)$$

where, $P_{\sigma(k)}$ are positive definite matrices $\forall\ \sigma(k) \in \Lambda$.

Let the value of $\sigma(k)$ at the $k^{th}$ and $(k+1)^{th}$ time instant be $i$ and $j$ respectively, where $i, j \in \Lambda$. The difference of the Lyapunov function between the two instants of time is given by:

$$\begin{aligned}\Delta V(k) &= V(k+1) - V(k) \\ &= \Gamma^T(k+1) P_j \Gamma(k+1) - \Gamma^T(k) P_i \Gamma(k) \\ &= \Gamma^T(k)\left(\Phi_i^T P_j \Phi_i - P_i\right)\Gamma(k)\end{aligned} \quad (7)$$

For any $\Gamma(k) \neq 0$, $\Delta V(k) < 0$ if (5) holds.

Thus, $\lim_{k \to \infty} V(k) = 0\ \forall\ k \in \mathbb{Z}_+$. Hence the system (3) is asymptotically stable. □

## IV. CONTROLLER DESIGN FOR ROBUST STABILITY

The controller design problem is recast into an optimization problem and is solved with the help of a stochastic optimization algorithm known as the Particle Swarm Optimization. The optimization variables are the predictive gains of the controller, viz. $K_l$ where $l \in \{2, 3, ..., M_{drop}\}$. The objective function to be maximized is as follows:

$$J = \underbrace{\iint \cdots \int f\left(K_{11}, \cdots, K_{1n}\right) dK_{11} \cdots dK_{1n}}_{n-fold} \quad (8)$$

where, $f$ is the region enclosed by the components of the nominal feedback controller gain $K_1$ satisfying equation (5). For a system consisting of 2 and 3 states (8) represents the area and volume enclosed by the components of $K_1$ respectively.

For the corresponding minimization problem, the negative value of the objective function is taken. Once the predictive gains are found out, the largest hyper-sphere is fitted in the stable region by numerical techniques and the co-ordinates of the centre give the nominal gains of the controller. $J$ is numerically found out by summing the number of points which satisfies (5), evaluated at equally spaced discrete points in the n-dimensional nominal controller space $K_1$.

### A. Particle Swarm Optimization Algorithm

The PSO tries to optimize an objective function $f(x)$ with respect to the design variable $x \in \mathbb{R}^n$. It is expressed as:

$$\underset{x \in \mathbb{R}^n}{\text{minimize}}\ f(x) \quad (9)$$

where, the objective function $f : \mathbb{R}^n \to \mathbb{R}$ and the n-dimensional search space $\mathbb{G} \in \mathbb{R}^n$ is pre-specified by the user. The PSO algorithm consists of a swarm of particles $x_i\ \forall\ i \in \{1, 2, ..., n_p\}$. The maximum number of particles

$n_p$ is specified by the user. The particles $x_i$ search for an optimal solution $x' \in \mathbb{R}^n$ of (9). The position of the $i^{th}$ particle is denoted by $x_i := (x_{i,1}, x_{i,2}, ..., x_{i,n})^T \in \mathbb{R}^n$ and the velocity is denoted by $v_i := (v_{i,1}, v_{i,2}, ..., v_{i,n})^T \in \mathbb{R}^n$, where $i \in \{1, 2, ..., n_p\}$. The position and velocity of the $i^{th}$ particle $x_i \in \mathbb{R}^n$ is updated in each iteration, based on the following equations, for $k \in \mathbb{Z}_+$ which indicates the iteration number.

$$x_i^{k+1} = x_i^k + v_i^{k+1} \quad (10)$$

$$v_i^{k+1} = \alpha v_i^k + \beta_1 \theta_{1,i}^k (x_i^{best,k} - x_i^k) + \beta_2 \theta_{2,i}^k (x_{swarm}^{best,k} - x_i^k) \quad (11)$$

where $\alpha$ is the inertia factor, $\beta_1$ is the cognitive learning rate and $\beta_2$ is the social learning rate and is pre-specified by the user. These influence the exploration and exploitation properties of the particles and must be properly chosen [19] for quicker convergence. $\theta_{1,i}^k$ and $\theta_{2,i}^k$ represent random numbers uniformly distributed in the interval $[0,1]$. $x_i^{best,k}$ in (11) refers to the previously obtained best position of the $i^{th}$ particle and $x_{swarm}^{best,k}$ denotes the best position of the swarm at the current iteration $k$. This is expressed as:

$$x_i^{best,k} := \arg\min_{x_i^j}\{f(x_i^j), 0 \le j \le k\} \quad (12)$$

$$x_{swarm}^{best,k} := \arg\min_{x_i^k}\{f(x_i^k), \forall i\} \quad (13)$$

The pseudo-code for the PSO algorithm can be summarized as follows:

[*Step 1*]: Initialize $n_p$ particles randomly distributed in the search space $\mathbb{G} \in \mathbb{R}^n$ and calculate the objective function values for each particle $x_i \; \forall \; i \in \{1, 2, ..., n_p\}$. Set, $k = 0$. Determine $x_i^{best,0}$ and $x_{swarm}^{best,0}$.

[*Step 2*]: If the criteria for termination is satisfied, the algorithm terminates with the solution $x' := \arg\min_{x_i^j}\{f(x_i^j), \forall i, j\}$. Otherwise go to [Step 3].

[*Step 3*]: Use (10) and (11) to update the position and velocity of the particles and evaluate the corresponding objective functions at each position. Set, $k = k+1$. Determine $x_i^{best,k}$ and $x_{swarm}^{best,k}$ and go to [Step 2].

The termination criterion is set as the user specified maximum number of iterations $k_{max}$. The simulations for the illustrative examples are run with a population of 20 and $k_{max} = 100$. The scaling factors are chosen to be $\beta_1 = \beta_2 = 1.49$. The inertia factor is linearly varied from $\alpha = 0.9$ to $\alpha = 0.1$. These settings are similar to that in [20] and are seen to work well for typical controller design tasks in NCS applications. Other small changes in parameters give similar results. Rather the convergence and run time of the algorithm is affected more.

## V. ILLUSTRATIVE EXAMPLE

### A. Simulation Study

Let us consider the practical networked DC motor system as in [5]. With the sampling time chosen as 0.05 seconds, the discrete time plant can be expressed as:

$$x_p(k+1) = \begin{bmatrix} 1 & 0.004601 \\ 4.601 \times 10^{-3} & 4.018 \times 10^{-5} \end{bmatrix} x_p(k) + \begin{bmatrix} 0.3487 \\ 7.681 \end{bmatrix} u(k) \quad (14)$$

where, $x_p = [\theta_p \; \omega_p]^T$ and $\theta_p$ is the angular position and $\omega_p$ is the angular speed of the motor. We assume that $M_{drop} = 3$, i.e. at the most two consecutive packets can be dropped. The stability region for this system for constant values of state feedback gains without any prediction is as shown below.

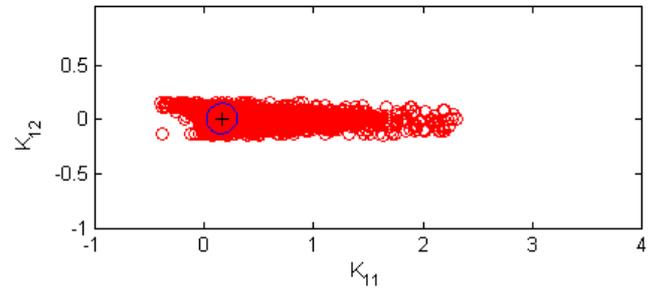

Fig. 3. Stability region with the largest circle fitted for simple state feedback controller.

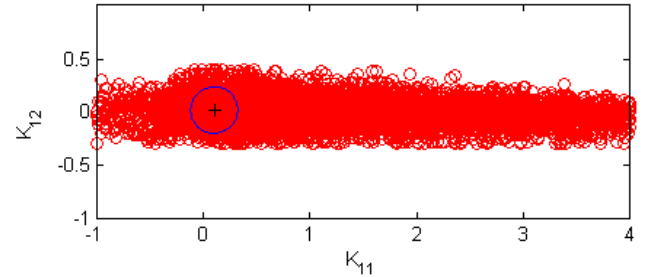

Fig. 4. Stability region with the largest circle fitted for variable gain state feedback controller.

The centre of the largest circle fitted in the maximum stabilizing region gives the value of gains as $K = [K_{11} \; K_{12}] = [0.16 \; 0.01]$. The predictive controller design methodology as exposited in Section IV yields the following values of the predictive state feedback gains $K_2, K_3$ which maximize the area of the nominal gain $K_1$:

$K_2 = [0.0837 \; 0.056]$ and $K_3 = [0.1527 \; 0.0682]$.

From the circle center the nominal gains are obtained as $K_1 = [K_{11}^{centre} \; K_{12}^{centre}] = [0.1 \; 0.02]$. The radius of the circle with constant gains in Fig. 3 is 0.14 and the radius of the circle with the predictive gains in Fig. 4 is 0.22. This shows that the predictive methodology allows for more variation in controller parameters while preserving stability in the sense of Lyapunov. The asymptotic stability of the system for initial condition $x_0 = [3.65\ 0]^T$ is shown in Fig. 5 and Fig. 6.

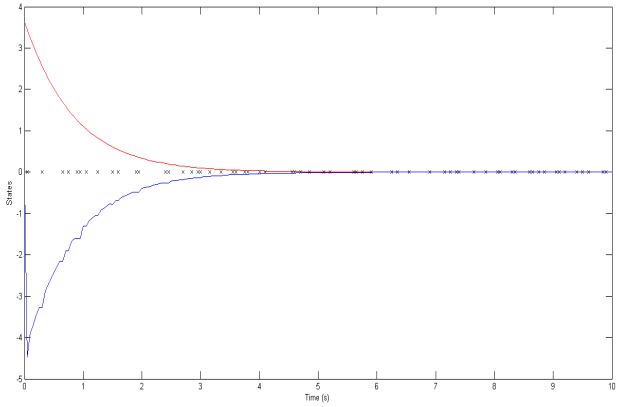

Fig. 5. Asymptotic stability for the system with simple state feedback controller, the crosses on the time axis represent packets which are dropped.

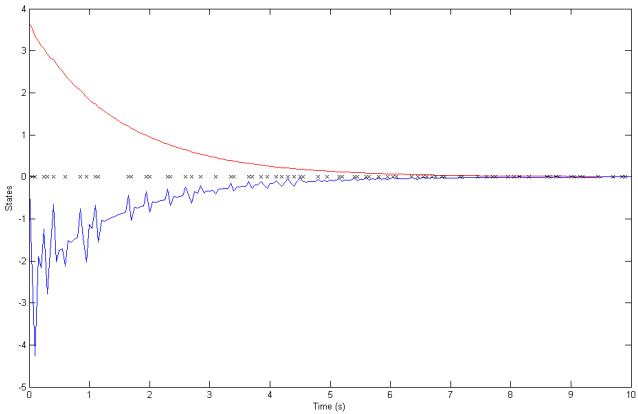

Fig. 6. Asymptotic stability for the system with predictive state feedback controller, the crosses on the time axis represent packets which are dropped.

*B. Few Discussions*

It is evident from Figs. 3 and 4, that the stability region in the sense of Lyapunov for the proposed predictive controller with different state feedback gains for different levels of packet dropouts is higher than the case with the simple state feedback controller. Since uncertainty in plant parameters is ultimately reflected in the change in controller parameters [16]-[17], this predictive approach gives more robustness to errors in plant parameter estimation. The actual gains are chosen to be the centre of the largest circle that can be fitted inside the stability region. This ensures maximum variation in controller parameters of both the states. The methodology can easily be extended to the case with higher number of state feedback gains. In that case the controller parameters can be found by fitting the largest hyper-sphere in n-dimensional parameter space. However pictorial representation in this case would be difficult.

Fig. 3 and 4 shows that the system can tolerate higher deviation in the state feedback gain $K_{11}$ associated with the position state ($\theta_p$) than that with $K_{12}$. For more realistic design it is possible to fit an ellipse instead of a circle in the stability region which permits more deviation in one state-feedback gain than the other. The obtained stabilizing gains are both positive and negative in Fig. 3 and 4. This implies that even a positive feedback action can stabilize a state feedback control loop since it acts as a switched system in the presence of packet losses. This is intrinsically different from the conventional linear time invariant (LTI) systems which can not be stabilized in such cases with a positive feedback in the sense of Lyapunov. Fig. 3 and 4 also shows that the stabilizing region is composed of discontinuous patches unlike [16] for LTI systems which is a typically due to the switched nature of the set of linear systems with packet drops.

Figs. 5 and 6 show the time domain asymptotic stability of the plant in the presence of packet dropouts for the case of constant feedback gains and predictive feedback gains respectively. The state responses in Fig. 6 are more jagged with sharp transitions which are not seen in Fig. 5. Thus it is evident that increase in stability region for the predictive controller comes at the cost of decrease in time domain performance.

VI. CONCLUSION

Predictive control methodology with constant state feedback gains are proposed in this paper. This takes the advantage of predictive control while reducing the computational complexity of calculation in each iteration, unlike conventional model predictive control approaches. The predictive gains are chosen with PSO to provide maximum parameter variation of the nominal controller, while ensuring stability in the Lyapunov sense in the presence of packet drops. The nominal gains are then chosen to be the centre of the largest possible circle in the stability region. Time domain simulation show that asymptotic stability is achieved for both cases but the increase in stability region for the predictive controller comes at the cost of decreased time domain performance, which is evident due to the presence of sharp transitions in the states. Nevertheless, this approach takes care of uncertainties in plant parameter variations and ensures robust stability. Scope of future work might include incorporating time domain performance while preserving robust stability conditions.

REFERENCES

[1] Joao P. Hespanha, Payam Naghshtabrizi, and Yonggang Xu, "A survey of recent results in networked control systems", *Proceedings of the IEEE*, vol. 95, no. 1, pp. 138-162, Jan. 2007.
[2] John Baillieul and Panos J. Antsaklis, "Control and communication challenges in networked real-time systems", *Proceedings of the IEEE*, vol. 95, no. 1, pp. 9-28, Jan. 2007.
[3] Wei Zhang, Michael S. Branicky, and Stephen M. Phillips, "Stability of networked control systems", *IEEE Control Systems Magazine*, vol. 21, no. 1, pp. 84-99, Feb. 2001.
[4] Junlin Xiong and James Lam, "Stabilization of linear systems over networks with bounded packet loss", *Automatica*, vol. 43, no. 1, pp. 80-87, Jan. 2007.
[5] Hongbo Li, Mo-Yuen Chow, and Zengqi Sun, "Optimal stabilizing gain selection for networked control systems with time delays and packet losses", *IEEE Transactions on Control Systems Technology*, vol. 17, no. 5, pp. 1154-1162, Sept. 2009.

[6] Junxia Mu, G.P. Liu, and David Rees, "Design of robust networked predictive control systems", *Proceedings of the 2005 American Control Conference*, vol. 1, pp. 638-643 June 2005.

[7] Senchun Chai, Guo-Ping Liu, David Rees, and Yuanqing Xia, "Design and practical implementation of internet-based predictive control of a servo system", *IEEE Transactions on Control Systems Technology*, vol. 16, no. 1, pp. 158-168, Jan. 2008.

[8] Yun-Bo Zhao, Guo-Ping Liu, and David Rees, "Networked predictive control systems based on the Hammerstein model", *IEEE Transactions on Circuits and Systems-II: Express Briefs*, vol. 55, no. 5, pp. 469-473, May 2008.

[9] Yun-Bo Zhao, Guo-Ping Liu, and David Rees, "A predictive control-based approach to networked Hammerstein systems: design and stability analysis", *IEEE Transactions on Systems, Man, and Cybernetics-Part B: Cybernetics*, vol. 38, no. 3, pp. 700-708, June 2008.

[10] Wenshan Hu, Guo-Ping Liu, and David Rees, "Event-driven networked predictive control", *IEEE Transactions on Industrial Electronics*, vol. 54, no. 3, pp. 1603-1613, June 2007.

[11] Guo-Ping Liu, Yuanqing Xia, David Rees, and Wenshan Hu, "Design and stability criteria of networked predictive control systems with random network delay in the feedback channel", *IEEE Transactions on Systems, Man, and Cybernetics-Part C: Applications and Reviews*, vol. 37, no. 2, pp.173-184, March 2007.

[12] Guo-Ping Liu, Yuanqing Xia, Jie Chen, David Rees, and Wenshan Hu, "Networked predictive control of systems with random network delays in both forward and feedback channels", *IEEE Transactions on Industrial Electronics*, vol. 54, no .3, pp. 1282-1297, June 2007.

[13] Yun-Bo Zhao, Guo-Ping Liu, and David Rees, "Design of a packet-based control framework for networked control systems", *IEEE Transactions on Control Systems Technology*, vol. 17, no. 4, pp. 859-865, July 2009.

[14] Hai Lin, Guisheng Zhai, and Panos J. Antsaklis, "Robust stability and disturbance attenuation analysis of networked control systems", *Proceedings of 42$^{nd}$ IEEE Conference on Decision and Control, 2003*, vol. 2, pp. 1182-1187, Dec. 2003.

[15] Marieke Cloosterman, Nathan van de Wouw, Maurice Heemels, and Henk Nijmeijer, "Robust stability of networked control systems with time-varying networked-induced delays", *2006 45$^{th}$ IEEE Conference on Decision and Control*, pp. 4980-4985, Dec. 2006, San Diego.

[16] Ming-Tzu Ho, Aniruddha Datta, and S.P. Bhattacharyya, "Robust and non fragile PID controller design", *International Journal of Robust and Nonlinear Control*, vol. 11, no. 7, pp. 681-708, June 2001.

[17] L. Ou, W. Zhang, and D. Gu, "Sets of stabilising PID controllers for second-order integrating processes with time delay", *IEE Proceedings - Control Theory and Applications*, vol. 153, no. 5, pp. 607-614, Sept. 2006.

[18] T.C. Yang, C. Peng, D. Yue, and M.R. Fei, "New study of controller design for networked control systems", *IET Control Theory & Applications*, vol. 4, no. 7, pp. 1109-1121, July 2010.

[19] James Kennedy and Russell Eberhart, "Particle swarm optimization", *Proceedings of IEEE International Conference on Neural Networks, 1995*, vol. 4, pp. 1942-1948, Nov.-Dec. 1995, Perth, WA, Australia.

[20] Indranil Pan, Saptarshi Das, and Amitava Gupta, "Tuning of an optimal fuzzy PID controller with stochastic algorithms for networked control systems with random time delay", *ISA Transactions*, vol. 50, no. 1, pp. 28-36, Jan. 2011.